\documentclass[11pt]{article}
             \newlength{\oldparindent}
             \usepackage{amsmath}
             \usepackage{latexsym}
             \usepackage{t1enc}
             \usepackage[latin1]{inputenc}
             \usepackage[french]{babel}
          
           \newtheorem{theoreme}{Th\'eor\`eme}[section]

             \newtheorem{definition}[theoreme]{Definition}

            \newtheorem{lemme}[theoreme]{Lemme}
             \newtheorem{remarque}[theoreme]{Remarque}
             
             \newtheorem{hyp}[theoreme]{Hypoth\`ese}

             \font\myf=msbm10 at 12pt
             \newcommand{\R}{\hbox{\myf R}}
             
             \newcommand{\pr}{\hbox{\myf P}}
             \newcommand{\esp}{\hbox{\myf E}}

             \newcommand{\F}{{\cal F}}

\begin{document}

        \def\second {\vtop{\baselineskip=11pt
             %\hbox to 125truept{\hss  \hss}}}
             \hbox to 125truept{\hss  Diana DOROBANTU\footnote{Laboratoire de Statistique et Probabilit\'es à l'Université Toulouse 3, dorobant@cict.fr} \hss}\vskip2truept}}

            \title{{ Arr\^et optimal pour les processus de Markov forts et les fonctions affines }}
    \author{ \second}
             \date{}
             \maketitle

 \vspace{1cm}
             \hbox to 125truept{\hss \bf  Résumé \hss}\vskip2truept
             
 Dans cette Note  nous \'etudions des probl\`emes d'arr\^et optimal pour des processus de Markov forts et des fonctions affines.  Nous donnons une justification de la forme de l'enveloppe de Snell en utilisant les r\'esultats classiques d'arr\^et optimal. Nous justifions \'egalement  la convexit\'e de la fonction valeur et sans nous restreindre a priori \`a une classe particuli\`ere de temps d'arr\^et, nous en d\'eduisons que le plus petit temps d'arr\^et  optimal est n\'ecessairement un temps d'atteinte. 
   \vskip 0.5\baselineskip

 \vspace{1cm}
             \hbox to 125truept{\hss \bf Abstract \hss}\vskip2truept
             
{\bf Optimal stopping for strong Markov processes and affine functions }

 \vspace{0.25cm}
 
In this Note we study optimal stopping problems for  strong Markov processes and  affine functions.  We give a justification of  the  Snell envelope form using  standard results of optimal stopping. We also justify the convexity  of the value function, and without  a priori restriction to a particular class of stopping times, we deduce that the   smallest  optimal stopping time is necessarily a hitting time. 

\section*{Abridged English version}
% Texte de la version abrégée en anglais
         
         We choose to solve a particular optimal stopping problem for strong Markov processes. Without a priori restriction  to a particular class of stopping times, we propose a method to find the optimal stopping time form (it will be a hitting time).

            In fact we seek to control a stochastic process V of the form $V=ve^X$ where $v$ is a real strictly positive constant and $X$ a  strong Markov process such that $X_0=0$. We consider the following optimal stopping problem :
$$J_t=esssup_{ \tau\in \Delta, \tau\geq t} \esp\left(e^{-r\tau}f(V_{\tau})\mid \F^V_t\right),$$
 where   $r>0$,  $\F^V_t=\sigma(V_s, ~s\leq t)$,  $\Delta$ is the set of   
$\F^V_.$-stopping times and $f$ is a decreasing affine function.  

 We suppose that the process $X$ checks the following  assumptions :
 
  \textbf{H1 :}  $\pr(lim_{t \downarrow 0} X_t=X_0)=1.$
 
 \textbf{H2 :} The process $(e^{-rt+X_t}, ~t\geq 0)$ is of class $D$.
 
  \textbf{H3 :}   $inf_{t\geq 0}e^{-rt}\esp(e^{X_t})=0$.
 
 \textbf{H4 :} The support of $X_t$ is $\R$ for all $t>0$.

 We seek a stopping time  $\tau^*$ which maximizes $\tau\mapsto \esp(Y_{\tau}\mid \F^V_t)$ where the process $Y$ has the form  \\$t\mapsto Y_t=e^{-rt}f(V_{t})$. In many papers the optimal stopping time is supposed from the beginning to be a hitting time, here we show that the optimal stopping time is necessarily of the form $inf\{t\geq 0 : ~V_t\leq b\}$. Under Assumptions $H1-H4$, we give a justification of the form of the Snell envelope of the process $Y$ using  standard results of optimal stopping of \cite{KLM} : $J$ has the form $t\mapsto J_t=e^{-rt}s(V_t)$ where the function $s$ is called "r-reduite" of $f$. We also argue the convexity of the function $s$. Using a result of  optimal stopping of \cite{Sh}, the smallest optimal stopping time  has the form    $$\tau^*=inf\{u : f(V_u)=s(V_u)\}.$$ The main result is given by Theorem \ref{th1}  which allows to deduce from the convexity of $s$, that  $\tau^*$ is necessarily a hitting time. 
\selectlanguage{francais}
% texte principale
\section{Introduction}
\label{}
  
               Nous nous pla\c{c}ons sur un espace de probabilit\'e filtr\'e $(\Omega,\F,(\F_t)_{t\geq 0},\pr)$ sur lequel  nous consid\'erons $V$ un processus stochastique s'\'ecrivant sous la forme $$V=ve^X$$ o\`u 
             $v$ est une constante r\'eelle strictement positive  et $X$ est un processus de Markov fort tel que $X_0=0$.
      Afin d'enlever l'ambigu\"it\'e,  nous utiliserons parfois la notation $V^v=ve^X$, pour tout $v>0$.
      
       Dans la suite  $\esp(.|V_0=v)$ et $\pr(.|V_0=v)$ sont not\'es $\esp_v(.)$ et $\pr_v(.)$.

               Nous introduisons $\F^V$ la filtration compl\'et\'ee c\`ad, engendr\'ee par le processus $V$, 
             $\F^V_t=\sigma(V_s, ~s\leq t)$ et nous consid\'erons le probl\`eme d'arr\^et optimal suivant :
\begin {equation}
             \label{2}
 J_t=esssup_{\tau\in \Delta, \tau\geq t} \esp\left(e^{-r\tau}f(V_{\tau})\mid \F^V_t\right),
             \end {equation}
             o\`u  $r>0$,  $\Delta$ est l'ensemble des  
             $\F^V_.$-temps d'arr\^et et $f$ est une fonction affine d\'ecroissante de la forme $f(v)=-\alpha v+c$ o\`u $\alpha>0$, $c>0$.
              \begin{definition}
             Le temps d'arr\^et $\tau^*_t$ est optimal s'il maximise (\ref{2}), c'est \`a dire 
             
             $$\esp\left(e^{-r\tau^*_t}f(V_{\tau^*_t})\mid \F^V_t\right)=esssup_{\tau\in \Delta, \tau\geq t}\esp\left(e^{-r\tau}f(V_{\tau})\mid \F^V_t\right).$$
             \end{definition}
             
             Nous supposons  que le processus $X$ v\'erifie les  hypoth\`eses  suivantes :
             
            \begin{hyp}
             \label{h0} 
             $\pr(lim_{t \downarrow 0} X_t=X_0)=1.$
             \end{hyp}
             
      \begin{hyp}
             \label{h1}
             Le processus $(e^{-rt+X_t}, ~t\geq 0)$ est de classe $D$.
             \end{hyp}
             
             \begin{hyp}
             \label{h1i}
            $inf_{t\geq 0}e^{-rt}\esp(e^{X_t})=0$.
             \end{hyp}
             
             \begin{hyp}
             \label{h2}
             Le support de la loi de $X_t$ est $\R$ pour tout $t>0$.
             \end{hyp}
             
            Dans la suite nous proposons une m\'ethode qui permet de trouver la forme du plus petit temps d'arr\^et optimal du probl\`eme (\ref{2}). Nous avons appliqu\'e cette m\'ethode \`a des processus de L\'evy \cite{d1,d2}, mais cela s'\'etend \`a des processus plus g\'en\'eraux.
            
                \section{Le temps d'arr\^et optimal}
            \label{section3.1_ch1}

             Dans cette partie nous montrons que le probl\`eme (\ref{2}) admet au moins un temps d'arr\^et optimal et que le plus petit temps d'arr\^et optimal est un temps d'atteinte.

         Sous l'hypoth\`ese \ref{h1}, le processus $(t\mapsto e^{-rt}f(V_t), ~t\geq 0)$ est de classe $D$. D'apr\`es le th\'eor\`eme   3.4 de   \cite{KLM}, l'enveloppe de Snell $J$ du processus
             $Y$  est de  la forme $\left(e^{-rt}s(V_t)\right)_{t\geq 0}$.  De plus, comme le support de   $V_t$ est $\R^*_+$  (c'est une cons\'equence directe de l'hypoth\`ese \ref{h2}), alors la d\'efinition  (\ref{2}) prise en $t=0$ donne  $J_0=s(v)$. Il suffit donc de regarder le probl\`eme en $t=0$ 
             
             \begin{equation}
             \label{3}
             s(v)=sup_{\tau\in \Delta} \esp_v\left(e^{-r\tau}f(V_{\tau})\right).
             \end{equation}

                    \begin{theoreme}
             \label{th1}
               Sous  les hypoth\`eses \ref{h0}, \ref{h1}, \ref{h1i} et \ref{h2}, il existe  au moins un temps d'arr\^et optimal pour le probl\`eme (\ref{3}).
                
                 Pour tout $c>0$, 
             il existe $b_c>0$  tel que le plus petit temps d'arr\^et optimal est de 
             la forme
             $$\tau_{b_c}=inf \{t\geq 0 : ~V_t \leq b_c\}.$$
             \end{theoreme}

             La d\'emonstration de ce th\'eor\`eme n\'ecessite plusieurs r\'esultats.
 
  La fonction $s$ est une fonction convexe (d\'ecroissante)  car c'est le sup de fonctions affines  (d\'ecroissantes) :
             $$s(v)=sup_{\tau\geq 0}\esp_v\left[e^{-r\tau}(-\alpha 
             V_{\tau}+c)\right]=sup_{\tau\geq 0}
             \esp_1\left[e^{-r\tau}(-\alpha  vV^1_{\tau}+c)\right].$$

             \begin{remarque}
             La fonction $s$ \'etant convexe, elle est donc continue.
             \end{remarque}

          Remarquons que $s$ est une fonction positive puisque
             $$s(v)\geq sup_{t\geq 0} \esp_v\left[e^{-rt} (-\alpha V_t+c)\right]\geq 
             sup_{t\geq 0} \esp_v\left[-e^{-rt}\alpha V_{t}\right]=sup_{t\geq 0}-\alpha v\esp\left[e^{-rt+X_t}\right]=0,$$
             o\`u pour la derni\`ere \'egalit\'e nous avons utilis\'e l'hypoth\`ese \ref{h1i}.

 Le th\'eor\`eme 3.3 page 127   de \cite{Sh}, permet de trouver le  temps d'arr\^et optimal d'un probl\`eme du type $sup_{\tau\geq 0}\esp_v\left[f(V_{\tau})\right]$ o\`u $f$ est une fonction mesurable. En g\'en\'eralisant ce r\'esultat, nous en d\'eduisons facilement que le th\'eor\`eme 3.3 page 127   de \cite{Sh} peut \^etre appliqu\'e \`a des processus de la forme $t\mapsto e^{-rt}f(V_t)$. Nous ne pouvons pas appliquer directement ce r\'esultat pour le probl\`eme (\ref{3}) puisque le processus  \\$t\mapsto e^{-rt}f(V_t)$ ne v\'erifie pas les hypoth\`eses du th\'eor\`eme ; par cons\'equent nous sommes oblig\'es de r\'e\'ecrire la fonction $s$ sous une autre forme. 
 
             \begin{lemme}
             \label{lem3}
             Soit pour $v>0$
             $$s(v)=sup_{\tau\geq 0}\esp_v\left[e^{-r\tau}(-\alpha 
             V_{\tau}+c)\right]~~\hbox{et}~~s^+(v)=sup_{\tau\geq 0}\esp_v\left[e^{-r\tau}(-\alpha V_{\tau}+c)^+\right],$$
             o\`u $x^+=max(x, 0)$.
Sous les hypoth\`eses \ref{h0}, \ref{h1}, \ref{h1i} et \ref{h2}   $s^+(v)>0$ et $s(v)=s^+(v)$ pour  tout $v>0$.
             \end{lemme}
             
             \begin{preuve}
Nous montrons que s'il existe $v_0>0$ tel que $s(v_0)<s^+(v_0)$, alors il existe $v_1>0$ tel que $s^+(v_1)=0$. Nous montrons ensuite, que cette derni\`ere relation ne peut pas \^etre satisfaite. 
             
            Par construction, pour tout $v>0$, $s(v)\leq s^+(v)$.
Supposons qu'il existe $v_0>0$ tel que $s(v_0)<s^+(v_0)$.
             
      Sous l'hypoth\`ese \ref{h0}, le processus $V_.$ est continu \`a droite en $0$, le processus $t \rightarrow Y^+_t=e^{-rt}(-\alpha V_t+c)^+$ est \`a  valeurs dans $[0, ~c]$, alors les hypoth\`eses du th\'eor\`eme 3.3 page 127   de \cite{Sh}  sont v\'erifi\'ees pour $Y^+$.
On note $f^+(v)=(-\alpha  v+c)^+$ ;  le temps d'arr\^et
              $$\tau^+=inf\{u\geq 0 : f^+(V^{v_0}_u)=s^+(V^{v_0}_u)\}$$ est le plus petit temps  d'arr\^et optimal du 
             probl\`eme $$s^+(v_0)=sup_{\tau\geq 0}\esp_{v_0}\left[e^{-r\tau}(-\alpha V_{\tau}+c)^+\right].$$

             D'apr\`es la d\'efinition de $s$ et $s^+$ :
             $$\esp_{v_0}\left[ e^{-r\tau^+}f(V_{\tau^+})\right]\leq s(v_0)<s^+(v_0)=\esp_{v_0}\left[ e^{-r\tau^+}f^+(V_{\tau^+})\right]$$
             et par suite $$\esp_{v_0}\left[ e^{-r\tau^+}\left(f(V_{\tau^+})-f^+(V_{\tau^+})\right)\right]<0, ~\pr_{v_0}\left(\{\omega : f(V_{\tau^+})<0\}\right)>0
            \hbox{~et~} \pr_{v_0}\left(\{\omega : s^+(V_{\tau^+})=0\}\right)>0.$$ 
             
           Donc il existe $v_1$ tel que $s^+(v_1)=0$.  Alors pour tout temps d'arr\^et $\tau$, $\pr_{v_1}$-presque s\^urement $e^{-r\tau}f^+(V_{\tau})=0$ et en particulier pour tout $t\in \R_+$, $f^+(V_t)=0$. Ceci entra\^ine que $\pr_{v_1}$-presque s\^urement $V_t\geq \frac{c}{\alpha}$ d'o\`u la contradiction puisque sous l'hypoth\`ese \ref{h2}, la loi de $V_t$ a pour support $\R^*_+$. Donc $s^+(v)>0$ pour tout $v\in\R^*_+$ et $s(v)=s^+(v)$.
             \end {preuve}
             \hfill  $\Box$
             \\
             
             D'apr\`es le lemme \ref{lem3}, le probl\`eme (\ref{3}) se ram\`ene \`a un probl\`eme d'arr\^et optimal pour une option Put am\'ericaine. Un tel probl\`eme a \'et\'e \'etudi\'e par plusieurs auteurs quand $X$ est un processus de L\'evy (voir par exemple \cite{boy,kyp,mo,pham}). Dans la suite,  nous utilisons un raisonnement proche de celui utilis\'e par Pham \cite{pham}. Il se place dans un cadre plus restreint que le n\^otre et \'etudie un probl\`eme d'arr\^et optimal \`a horizon fini pour une option Put am\'ericaine quand  $X$ est un processus de L\'evy. 
            \\
          
   \begin{preuve} (du th\'eor\`eme \ref{th1})
             \\
                    D'apr\`es le lemme \ref{lem3}, le probl\`eme (\ref{3}) peut se 
             r\'e\'ecrire sous la forme $sup_{\tau\geq 0}\esp(Y_{\tau }^+)$.  Les hypoth\`eses 
             du th\'eor\`eme 3.3 page 127   de \cite{Sh} sont
               v\'erifi\'ees et le temps d'arr\^et
               $$\tau^*=inf\{u\geq 0 :
              f^+(V_u)=s^+(V_u)\}$$ est le plus petit temps  d'arr\^et optimal.
Or $s(v)=s^+(v)>0$ pour tout $v>0$, d'o\`u 
              $$\tau^*=inf\{u\geq 0 :
              f(V_u)=s(V_u)\}$$ est le plus petit temps  d'arr\^et optimal. 
                 La fonction  $s$ est
               major\'ee par $c$ puisque  $Y_.^+$ est major\'e par
               $c$ et $lim_{v\downarrow 0} s(v)=lim_{v\downarrow 0} 
             f(v)=c$.
             
             La fonction $s$ \'etant  convexe et $f$ une fonction
                affine, alors
             $inf \{v>0 : f(v)<s(v)\}$ est \'egal \`a \\$sup \{v>0 : f(v)=s(v)\}$,  que 
             l'on note $b_c$. En effet, soit $b_c'=sup \{v : f(v)=s(v)\}$  et \\$b_c=inf \{v : f(v)<s(v)\}$. Comme $lim_{v\downarrow 0} s(v)=lim_{v\downarrow 0}  f(v)$, alors $b_c'$ existe et $b_c'\geq 0$. Si $b_c=0$, alors $b_c'=0.$

             Si  $b_c>0$, alors pour tout $v<b_c$, $f(v)=s(v)$ ; en particulier $f(b_c-\frac{1}{n})=s(b_c-\frac{1}{n})$. En faisant tendre $n$ vers l'infini, comme $s$ et $f$ sont continues, alors $f(b_c)=s(b_c)$, d'o\`u $b_c\leq b_c'$.  Supposons par l'absurde que $b_c<b_c'$, alors il existe $v$, $b_c<v<b_c'$ tel que $f(v)<s(v)$.  Or $s$ est convexe, donc d'apr\`es le lemme des trois cordes  :
             $$\frac{s(v)-s(b_c)}{v-b_c}\leq \frac{s(b_c')-s(v)}{b_c'-v}.$$
            Comme, par continuit\'e, $f(b_c')=s(b_c')$, alors
            $$\frac{s(v)-f(b_c)}{v-b_c}\leq \frac{f(b_c')-s(v)}{b_c'-v}.$$
            Comme $s(v)>f(v)$, alors
             $$\frac{f(v)-f(b_c)}{v-b_c}<\frac{s(v)-f(b_c)}{v-b_c}\leq \frac{f(b_c')-s(v)}{b_c'-v}<\frac{f(b_c')-f(v)}{b_c'-v},$$
 d'o\`u la contradiction puisque, comme $f$ est affine, alors  $\frac{f(v)-f(b_c)}{v-b_c}=\frac{f(b_c')-f(v)}{b_c'-v}=-\alpha$. Donc $b_c=b_c'$.
             
            Ceci signifie que le plus petit  temps d'arr\^et optimal
             $\tau^*$ est aussi le premier temps de passage dans l'intervalle $]0, ~b_c]$.
             \end {preuve}
             \hfill  $\Box$
             \\
             
             \begin{remarque}
             Le plus petit temps d'arr\^et optimal du probl\`eme (\ref{2}) est de la forme   $$inf \{u\geq t : ~V_u \leq b_c\}.$$
             \end{remarque}
       
       Soit $\tau_x=inf \{u\geq t : ~X_u \leq x\}.$ Lorsque les transform\'ees de  Laplace de $\tau_x$ et de  $(\tau_x, X_{\tau_x})$ sont connues, le seuil optimal peut \^etre calcul\'e explicitement \cite{d1,d2}.   
        
% etc, etc

% Les remerciements sont dans une section, sans numérotation

%\section*{Remerciements}
% Remerciements - texte ici

\end{document}